\documentclass[12pt]{article}
\usepackage{amsmath,amssymb,latexsym,a4wide,graphicx}

\newtheorem{theorem}{Theorem}
\newtheorem{proposition}[theorem]{Proposition}
\newtheorem{lemma}[theorem]{Lemma}
\newtheorem{corollary}[theorem]{Corollary}

\begin{document}

\title{\vspace{-2.3cm} Regular self-dual and self-Petrie-dual maps \\ of arbitrary valency}

\author{}
\date{}
\maketitle

\begin{center}
\vspace{-1.3cm}

{\large Jay Fraser, Olivia Jeans} \\
\vspace{1.5mm} {\small Open University, Milton Keynes, U.K.}\\

\vspace{5mm}

{\large Jozef \v Sir\'a\v n} \\
\vspace{1.5mm} {\small
Open University, Milton Keynes, U.K., and \\ Slovak University of Technology, Bratislava, Slovakia}

\vspace{4mm}

\end{center}

\begin{abstract}
The main result of D. Archdeacon, M. Conder and J. \v{S}ir\'a\v{n} [Trans. Amer. Math. Soc. 366 (2014) 8, 4491-4512] implies existence of a regular, self-dual and self-Petrie dual map of any given even valency. In this paper we extend this result to any odd valency $\ge 5$. This is done by algebraic number theory and maps defined on the groups ${\rm PSL}(2,p)$ in the case of odd prime valency $\ge 5$ and valency $9$, and by coverings for the remaining odd valencies.

\vskip 3mm


\end{abstract}

\vskip 3mm

\section{Introduction}

In this paper we consider regular maps (that is, cellular embeddings of graphs on closed surfaces) with the highest `level of symmetry', which are, in addition, invariant under the operators of duality and Petrie duality. Regular maps have been addressed in a number of papers and we refer here to the latest survey \cite{Si-surv} for a large number of details; here we just sum up the essentials needed for our purposes.
\smallskip

From an algebraic point of view, a regular map $M$ can be identified with a finite group $G$ with three distinguished involutory generators $x,y,z$ and relators $(yz)^k$, $(zx)^\ell$ and $(xy)^2$ so that $x$ and $y$ commute; we will formally write $M=(G;x,y,z)$ to encapsulate the situation. The pair $(k,\ell)$ is the {\em type} of $M$, and we will assume throughout that $k,\ell\ge 3$; the type is {\em hyperbolic} if $1/k+1/\ell < 1/2$. Geometrically and topologically, elements of $G$ may be identified with flags (which correspond to mutually incident vertex-edge-face triples) and left cosets of the subgroups $\langle x,y \rangle$, $\langle y,z\rangle$ and $\langle z,x\rangle$ represent edges, vertices and faces of the embedded graph, with incidence given by non-empty intersection of cosets. Moreover, left multiplication by elements of $G$ on the cosets induce map automorphisms of $M$ and, in fact, $G$ is isomorphic to the (full) automorphism group ${\rm Aut}(M)$ of $M$. Conjugates of $x$, $y$ and $z$, respectively, induce automorphisms that locally act on $M$ as reflections along some edge, in some edge, and in an axis of some corner of $M$. Similarly, conjugates of $r=yz$ and $s=zx$ represent rotations about vertices and face centres of the map; in particular, every vertex has valency $k$ and every face is bounded by a closed walk of length $\ell$. The map $M$ is {\em orientable} (meaning that its underlying surface is) if and only if $G^+=\langle r,s\rangle$ is a subgroup of $G$ of index two, and {\em non-orientable} otherwise. Thus, in the non-orientable case, the entire group $G$ can be generated by the two rotations $r$ and $s$ only, and the involutions $x,y,z$ are then expressible in terms of $r$ and $s$; in such a situation we also write $M=(G;r,s)$.
\smallskip

Since every automorphism of a map (thought of as a permutation of flags that preserves incidence along and across edges and within corners) is completely determined by its action on a single flag, in the `most symmetric' case when the automorphism group is transitive (and hence regular) on flags we may identify the group with the flag set, arriving at the description outlined above. But even in this `most symmetric' case the map may still admit `more symmetries'. Here we will consider two kinds of such `external symmetries', namely, invariance under the operators of duality and Petrie-duality of maps. The two operators are well known; informally, duality interchanges the roles of vertices and faces, and the Petrie dual of a map is formed by re-embedding its underlying graph so that the new faces are the left-right (`zig-zag') walks read out of the original map. A map is {\em self-dual} and {\em self-Petrie-dual} if it is isomorphic to its dual and Petrie dual, respectively. In the case of a {\em regular} map $M=(G;x,y,z)$ as above, it is also well known (cf. \cite{Si-surv}) that $M$ is self-dual if and only if the group $G$ admits an automorphism interchanging $x$ with $y$ and fixing $z$, and $M$ is self-Petrie-dual if $G$ has an automorphism interchanging $x$ with $xy$ and fixing $y$ and $z$. In \cite {ACS}, regular maps that are both self-dual and self-Petrie dual have been said to have {\em trinity symmetry}.
\smallskip

The natural question of the existence of regular maps of any valency with trinity symmetry was raised more than four decades ago. In \cite{W} it was suggested that the map $M=(G;x,y,z)$ for the group $G = \langle x,y,z; x^2,y^2,z^2, (xy)^2,(yz)^{2n},(zx)^{2n},(xyz)^{2n},(xzyzxyz)^2\rangle$ is a regular map with trinity symmetry, of valency $2n$ for every $n\ge 1$. This was eventually proved in \cite{ACS} in a much more general form, including also invariance under the so-called hole operators that represent additional levels of `external symmetries' not discussed here. {\sl However, the question remained almost completely open for odd-valent regular maps with trinity symmetry}. Note that such a map must necessarily be non-orientable because of self-Petrie-duality with Petrie walks of odd length. There is no such map of valency $3$ since the only regular map of type $(3,3)$ is the 2-skeleton of a tetrahedron, and the known examples of odd valency $k\ge 5$ were scarce and came out of three sources. Existence of regular maps with trinity symmetry of odd valency $k$ such that $5\le k\le 19$ was computationally verified in \cite{Con}. A method based on products of finite simple groups that could potentially lead to such maps for infinitely many odd values of $k$ was developed in \cite{Jon}, with explicit examples for $k=15$ (found also in \cite{ACS} by a different method) and $k=455$.
\smallskip

Here we completely settle the problem by showing that for every odd $k\ge 5$ there exists a regular, self-dual and self-Petrie-dual map of valency $k$. Our strategy is to establish this result first for every {\sl prime} $k\ge 5$ and for $k=9$ by algebraic methods motivated by those used in \cite{JMS} and applied to more detailed results of \cite{EHJ} on regular maps defined on linear fractional groups. We then extend this to non-prime odd values of $k\ge 5$ by an analogue of a covering tool from \cite{ACS}. The paper is organized accordingly: in sections \ref{sec:lfg} and \ref{sec:alg} we review results on regular self-dual and self-Petrie-dual maps on linear fractional groups and develop the algebraic methods needed for our purposes, and in section \ref{sec:main} we prove our general result and make a few concluding remarks.

\section{Regular maps on linear fractional groups}\label{sec:lfg}

Classification of all orientably-regular maps with orientation-preserving automorphism group isomorphic to ${\rm PSL}(2,q)$ or ${\rm PGL}(2,q)$ follows from \cite{Mac} and can be found in a somewhat more explicit form in  \cite{Sah}; the latter was re-interpreted and extended to regular maps (orientable or not) in \cite{CPS2}. Since we will be interested only in the special case of odd valency and face length we just reproduce the corresponding part of the classification result here (the cases when one of the entries in the type of the map is even are more involved and we refer to \cite{JMS} for details).

\begin{proposition}\label{prop:psl} Let $(k,\ell)\ne (5,5)$ be a hyperbolic pair with both entries odd and let $p$ be an odd prime dividing neither $k$ nor $\ell$. Let $e=e(k,\ell)$ be the smallest positive integer $j$ such that $2n\, | \, (p^j-\varepsilon_n)$ for each $n\in \{k,\ell\}$ and some $\varepsilon_n\in \{+1,-1\}$, and let $\xi_{n}$ be a primitive $2n$-th root of unity in ${\rm GF}(p^e)$ if $\varepsilon_n=1$ or in ${\rm GF}(p^{2e})$ if $\varepsilon_n=-1$. Further, let $D=\xi_k^2+\xi_k^{-2}+\xi_\ell^2+\xi_\ell^{-2}\ne 0$ and let
\[ R = \pm\left[\begin{array}{cc}  \xi_k & 0\\ 0 & \xi_k^{-1}
\end{array}\right] \ \ \ {\it and} \ \ \  S
=\pm (\xi_k-\xi_k^{-1})^{-1}\left[\begin{array}{cc}
-(\xi_\ell+\xi_\ell^{-1})\xi_k^{-1}\! & -D\\
1 & \!(\xi_\ell+\xi_\ell^{-1})\xi_k
\end{array}\right] \]
be elements of ${\rm PSL}(2,p^e)$ if $\varepsilon_k=1$ and of ${\rm PSL}(2,p^{2e})$ otherwise. Then,
\begin{itemize}
\item[{\rm (a)}] the group $G_{k,\ell}=\langle R,S\rangle$ is isomorphic to $PSL(2,p^{e})$, with $R$ of order $k$ and $S$ of order $\ell$;
\item[{\rm (c)}] $M=(G_{k,\ell};R,S)$ is a regular map of type $(k,\ell)$, which is non-orientable if and only if $-D$ is a square in $GF(p^e)$.
\end{itemize}
\end{proposition}

We note that if $p^e\equiv \pm 1$ (mod $10$), the group ${\rm PSL}(2,p^e)$ contains (up to conjugacy) two exceptional pairs $R,S$ as above for $(k,\ell)=(5,5)$ with the property that $\langle R,S\rangle \cong A_5$; this case (omitted from \cite[Theorem 2.2]{JMS}) is addressed in \cite{EHJ}. However, this situation does not apply in what follows.
\smallskip

Necessary and sufficient conditions for self-duality and self-Petrie-duality of the maps $M=(G_{k,\ell};R,S)$ from Proposition \ref{prop:psl} were established in \cite{EHJ}. As they are also quite complex we present here only a simple sufficient condition appearing as Corollary 4.3 in \cite{EHJ} which (in terms and notation of Proposition \ref{prop:psl}) can be re-stated as follows. 

\begin{proposition}\label{prop:ehj}
Let $k\ge 5$ be odd, and let $p\ge 5$ be a prime not dividing $k$. Further, let $\ell=k$ and let $\xi=\xi_k=\xi_\ell$ be a primitive $2k$-th root of unity in ${\rm GF}(p^e)$ or in ${\rm GF}(p^{2e})$ for $e=e(k,\ell)$ such that $3(\xi^2 + \xi^{-2})+2=0$. Then, $M=(G_{k,\ell};R,S)$ is a (non-orientable) self-dual and self-Petrie-dual regular map of valency $k$, with automorphism group isomorphic to ${\rm PSL} (2,p^e)$.
\end{proposition}

The condition $3(\xi^2 + \xi^{-2})+2=0$ is equivalent to $3(\xi+\xi^{-1})^2=4$ and for its fulfilment it is necessary that $3$ be a square in ${\rm GF}(p^e)$, $p\ge 5$. This, for $e=1$, holds if and only if $p\equiv \pm 1$ (mod $12$), and is always the case if $e\ge 2$. But we can say more. Namely, the element $\zeta = \xi^2$ in Proposition \ref{prop:ehj} is a primitive $k$-th root of unity in $F={\rm GF}(p^e)$ or $F={\rm GF}(p^{2e})$, and the condition $3(\zeta+\zeta^{-1})+2=0$ represents a quadratic equation in the prime field $F_p$ of $F$; it also says that $\zeta+\zeta^{-1}\in F_p$. The last fact is equivalent to  $(\zeta+\zeta^{-1})^p=\zeta+\zeta^{-1}$, which reduces to $(\zeta^{p-1}-1)(\zeta^{p+1}-1)=0$ in $F$. It follows that  either $\zeta\in F_p$ and $p\equiv 1$ (mod $2k$), or $\zeta$ lies in a quadratic extension of $F_p$ and $p\equiv -1$ (mod $2k$), and in both cases we have $e=1$ (recall that $k$ is assumed to be odd). The bulk of Proposition \ref{prop:ehj} may now be restated in a form more suitable for our future use.

\begin{corollary}\label{cor:ejh}
Let $k\ge 5$ be odd. Assume that there exists a prime $p\ge 5$ such that $p\equiv \pm 1\ ({\rm mod}\ 2k\ {\rm and}\ 12)$, and a primitive $k$-th root of unity $\zeta$ in a finite field of order $p$ or $p^2$ with the property $3(\zeta+\zeta^{-1}) +2=0$. Then, there exists a non-orientable self-dual and self-Petrie-dual regular map of valency $k$ with automorphism group ${\rm PSL}(2,p)$.
\end{corollary}

\section{Algebraic preliminaries}\label{sec:alg}

For any $k\ge 3$ let $\alpha$ be a primitive complex $k$-th root of unity; its minimal polynomial is the $k$-th cyclotomic polynomial. Let $h=\alpha+\alpha^{-1}$ and let $K={\mathbb Q}(h)$ be the field obtained by adjoining $h$ to the rationals. It is known \cite[Proposition 2.16]{Wash} that the ring ${\cal O}$ of algebraic integers of $K$ is ${\mathbb Z}(h)$. We will focus on the algebraic integer $g=3h+2\in {\cal O}$. Observe that $g\ne 0$, for otherwise $\alpha$ would be a root of a quadratic polynomial over ${\mathbb Z}$, contrary to $k\ge 3$.
\smallskip

Recall that the norm $N(y)$ of an element $y\in {\cal O}$ is defined as the product $\prod_t\sigma_t(y)$, where $\sigma_t$ denotes the injective homomorphism ${\cal O}\to {\mathbb C}$ into the field of complex numbers, uniquely determined by $\sigma_t(\alpha)=\alpha^t$, and $t$ ranges over all integers between $1$ and $(k-1)/2$ that are relatively prime to $k$. It is well known that $N(y)$ is an integer for any $y\in {\cal O}$, which is a consequence of the invariance of $N(y)$ under the endomorphisms $\sigma_t$.
\smallskip

For the norm of our element $g\in {\cal O}$ we thus have $N(g) = \prod_{t} (3\sigma_t(h)+2)$, the product being taken over all $t$ between $1$ and $(k-1)/2$, coprime to $k$. The $\varphi(k)/2$ images $\sigma_t(h)$ appearing in this product are precisely the roots of the minimal polynomial $\Psi(x)$ of degree $\varphi(k)/2$ for $h=\alpha+\alpha^{-1}$, see e.g. \cite{JMS}. So, if $\Psi(x) = \prod_t (x-\sigma_t(h)) = \sum_j a_jx^j$ where $j$ ranges from $0$ to $\varphi(k)/2$, then the integral coefficients $a_j$ will also appear in the expansion of the above product. More precisely, letting $r=\varphi(k)/2$ and $u=-2/3$, we have
\begin{equation}\label{eq:N1}
N(g) = \prod_{t} (3\sigma_t(h)+2) = (-3)^r\prod_t(u-\sigma_t(h)) = (-3)^r\sum_{j=0}^r a_ju^j = \sum_{j=0}^r (-3)^{r-j}2^{j}a_j \ .
\end{equation}

Let us consider what happens when we look at (\ref{eq:N1}) modulo $9$. Up to the last two terms all the remaining ones are a multiple of $9$ and so, noting that $a_r=1$,
\begin{equation}\label{eq:N2}
N(g) \equiv 2^r - 3{\cdot}2^{r-1}a_{r-1} \ \ ({\rm mod}\ 9)\ .
\end{equation}

We will show that if $k\ge 5$ and $k$ is a {\em prime}, then the norm $N(g)$ is not equal to $\pm 1$, which means that $g$ is then {\em not} a unit of the ring ${\cal O}$. Indeed, let $k\ge 5$ be a prime, so that   $r= \varphi(k)/2=(k-1)/2$. By \cite{SC} we then also have $a_{r-1}=1$, and the congruence (\ref{eq:N2}) becomes
\[ N(g) \equiv 2^{(k-1)/2} - 3{\cdot}2^{(k-3)/2} \equiv  -2^{(k-3)/2} \ \ ({\rm mod}\ 9)\ .\]
It is easy to check that $2^j\equiv \pm 1$ (mod $9$) for a positive integer $j$ if and only if $j$ is a multiple of $3$. This means that if $k$ is prime, the norm $N(g)$ can be congruent to $\pm 1$ (mod $9$) only if $(k-3)/2$ is a multiple of $3$, giving a contradiction if $k\ge 5$. Further, from (\ref{eq:N1}) with the help of $a_r=1$ and $a_0=\pm 1$ \cite{SC} it follows that if $k\ge 5$ is a prime, then $N(g)$ is divisible neither by $2$ nor by $3$. We thus have:

\begin{lemma}\label{lem:Ng} If $k\ge 5$ is a prime, then $N(g)\ne \pm 1$; in particular, the non-zero element $g\in {\cal O}$ is not a unit of the ring ${\cal O}$. Moreover, for every prime factor $p$ of $N(g)$ one has $p\ge 5$. \hfill $\Box$
\end{lemma}

Consider now the field $K'={\mathbb Q}(\alpha)$, an extension of $K$ of degree two. Let ${\cal O}'$ be the ring of algebraic integers of $K'$; it is well known \cite[Theorem 2.6]{Wash} that ${\cal O}'={\mathbb Z}(\alpha)$, and, of course, $[{\cal O}':{\cal O}]=2$. The (integral) norm $N'(z)$ of any $z\in {\cal O}'$ is now the product $\prod_t\sigma_t(z)$ taken over all injective homomorphism $\sigma_t:\ {\cal O}'\to {\mathbb C}$ given by $\sigma_t(\alpha)=\alpha^t$ for $t$ between $1$ and $k-1$ coprime to $k$, and one again has $N'(z)\in {\mathbb Z}$. The two norms, $N$ on ${\cal O}$ and $N'$ on ${\cal O}'$, are related by $N'(y) = (N(y))^2$ for each $y\in {\cal O}$.
\smallskip

We will keep assuming that $k\ge 5$ is an odd prime and we let $p\ge 5$ be an arbitrary prime divisor of $N(g)$, which exists by Lemma \ref{lem:Ng}. We continue by considering the ideal $\langle g,p\rangle$ of ${\cal O}'={\mathbb Z}(\alpha)$ generated by the elements $g$ and $p$.
\medskip

\begin{lemma}\label{lem:kp} If $k\ge 5$ is a prime and if $p\ge 5$ is a prime divisor of $N(g)$, the ideal $\langle g,p\rangle$ is proper in the ring ${\cal O}'$.
\end{lemma}

{\sl Proof.} Suppose that $\langle g,p\rangle = {\cal O}'$, which means that $1=Ag+ Bp$ for some $A,B\in {\cal O}'$. Clearly $A\ne 0$, for otherwise $1=N'(B)N'(p) = N'(B)p^{k-1}$ and so $N'(B)$ would not be an integer. Now, $1=N'(1)=N'(Ag+Bp) = \prod_{\sigma} \sigma(Ag+Bp)$, where the product is being taken over all the $\varphi(k)=k-1$ embeddings $\sigma:\ {\cal O}'\to {\mathbb C}$. Expansion of this product gives $N'(Ag+Bp) = N'(A)N'(g) + cp$ for some $c\in {\cal O}'$. Thus, $cp\in {\mathbb Z}$ and so either $c\in {\mathbb Z}$ or $c=\pm 1/p$. As $p$ is a divisor of $N'(g)=(N(g))^2$ and $N'(A)$ is a non-zero integer, in either case it follows that $N'(Ag+Bp)\ne 1$, a contradiction. \hfill $\Box$
\medskip

By Lemma \ref{lem:kp}, the ideal $\langle g,p\rangle$ is contained in some maximal ideal $J=J_p$ of the ring ${\cal O}'$. Since ${\cal O}'$ is a Dedekind domain, the ideal $J$ has finite index in ${\cal O}'$ and so ${\cal O}'/J$ is a finite field $F$ of characteristic $p$, that is, $F\cong GF(p^m)$ for some $m\ge 1$. Recalling our assumption of primality of $k$ we show that the (multiplicative) order of the element $\overline \alpha = \alpha + J$ in the field $F={\cal O}'/J$ is equal to $k$. Indeed, suppose this is not the case. Then, because of primality of $k$, the order of $\overline\alpha$ in $F$ would have to be one, meaning that $\overline\alpha=1$ in $F$. But then, since the element $\overline g = g + J$ is equal to zero in $F$, we would have $0 = \overline g = 3(\overline \alpha+ \overline \alpha^{\ -1}) + 2 = 8$ in $F$, a contradiction as $p$ is odd. Observe also that $k\ne p$ since no element in $F$ has multiplicative order $p$.
\smallskip

This way we have constructed a finite field $F$ of characteristic $p$ containing a primitive $k$-th root $\overline\alpha$ of unity such that $3(\overline \alpha + \overline {\alpha}^{\ -1}) + 2 = 0$. We now invoke
the analysis immediately preceding Corollary \ref{cor:ejh} in section \ref{sec:lfg}, which fully applies to our situation. As the result we conclude that $F$ is the prime field $F_p$ if and only if $\overline\alpha\in F_p$ for $p\equiv 1$ (mod $2k$); otherwise $F$ is a quadratic extension of $F_p$ for $p\equiv -1$ (mod $2k$). In both cases we have $p\equiv \pm 1$ (mod $12$) because $3$ has to be a square in $F_p$. Summing up, we have proved:
\medskip

\begin{proposition}\label{prop:pNg}
Let $k\ge 5$ be an odd prime and let $\alpha$ be a primitive complex $k$-th root of unity. Further, let $g=3(\alpha + \alpha^{-1})+2$ and let $N(g)$ be the norm of $g$ in the ring ${\mathbb Z}(\alpha)$. Then, $N(g)\notin \{0,\pm 1\}$, every prime divisor $p$ of $N(g)$ satisfies $p\ge 5$, $p\equiv \pm 1\ ({\rm mod}\ 2k\ {\rm and}\ 12)$, and for every such $p$ there is a finite field $F$ of order $p$ or $p^2$ containing a primitive $k$-th root $\overline \alpha$ of $1$ such that $\overline g = 3(\overline \alpha + \overline {\alpha}^{\ -1}) + 2 = 0$ in $F$. \hfill $\Box$
\end{proposition}

\section{The main result}\label{sec:main}

To obtain a restricted version of our main result for prime valencies at least five we just need to put the pieces together. Indeed, taking $\zeta=\overline\alpha$ in Proposition \ref{prop:pNg} and combining it with Proposition \ref{prop:ehj} and Corollary \ref{cor:ejh} immediately gives:

\begin{theorem}\label{thm:m1}
For every odd prime $k\ge 5$ there exists a prime $p\equiv \pm 1\ ({\rm mod}\ 2k\ {\rm and}\ 12)$ such that ${\rm PSL}(2,p)$ is the automorphism group of a (non-orientable) regular, self-dual and self-Petrie-dual map of valency $k$. \hfill $\Box$
\end{theorem}

We know that there is no $3$-valent regular map with trinity symmetry, but there is one of valency $3^2$ that can be constructed by the machinery of section \ref{sec:lfg} as follows. The element $2$ is a primitive $9$-th root of unity mod $73$, and so is $\zeta=2^4$ and its multiplicative inverse $\zeta^{-1}=2^5$, with $\zeta$ and $\zeta^{1}$ satisfying the condition $3(\zeta+\zeta^{-1})+2=0$ (mod $73$). By Proposition \ref{prop:ehj} the group ${\rm PSL}(2,73)$ carries a self-dual and self-Petrie-dual regular map of valency $9$.
\smallskip

Based on Theorem \ref{thm:m1} and the above remark we are now in position to prove a full version of our main result. As alluded to in the Introduction, this will be done with the help of coverings, and more specifically using a non-orientable analogue of Theorem 2.1 of \cite{ACS}. We state it here in a restricted version sufficient for our purpose.

\begin{theorem}\label{thm:acs}
If there is a non-orientable regular map of odd valency $d\ge 5$ with trinity symmetry and with automorphism group $G$, then for any odd integer $n\ge 3$ there is a non-orientable regular map of degree $nd$ with trinity symmetry and automorphism group isomorphic to $(\mathbb{Z}_n)^{1+|G|/4}\rtimes G$.
\end{theorem}

{\bf Proof (sketch).} As indicated, this result was proved in \cite[Theorem 2.1]{ACS} for orientable maps (and, in this category, in a much more general setting that included also external symmetries induced by hole operators). The parts of the proof in \cite{ACS} that refer to regularity, self-duality and self-Petrie-duality apply almost word-by-word to the non-orientable case and we thus give only a sketch of the arguments here. We will assume familiarity with the theory of lifts of maps by corner voltage assignments as explained e.g. in \cite{AGS,ARSS,ACS}; a {\em corner} of a regular map $M=(G;x,y,z)$ is any $2$-subset of the form $\{g,gz\}$ for $g\in G$.
\smallskip

Let now $M=(G;x,y,z)$ be a regular map as in the statement. For an odd $n\ge 3$ let $H={\mathbb Z}_n^{|G|/2}$ be the space of all $|G|/2$-tuples with entries from ${\mathbb Z}_n$ and let ${\cal E}$ be the set of unit vectors (those with exactly one non-zero coordinate, equal to $1$) in $H$. Define a corner voltage assignment $\sigma$ on flags of $M$ -- that is, on the elements of $G$ -- in the group $H$ by assigning the $|G|/2$ two-element subsets $\{\varepsilon, -\varepsilon\}$ for $\varepsilon \in {\cal E}$ to the $|G|/2$ corners $\{g,gz\}$ for $g\in G$ in an arbitrary one-to-one fashion. By arguments in the proof of Theorem 2.1 in \cite{ACS} that do not depend on orientability, the lift of the map $M$ of type $(d,d)$ by the voltage assignment $\sigma$ has $n^{-1+|G|/4}$ components, each isomorphic to a regular map $M^{\sigma} = (G^{\sigma}; x^{\sigma}, y^{\sigma}, z^{\sigma})$ of type $(nd,nd)$ for the group $G^{\sigma} = (\mathbb{Z}_n)^{1+|G|/4} \rtimes G$ and suitable involutory generators $x^{\sigma}, y^{\sigma}, z^{\sigma}$ of $G^{\sigma}$. Moreover, by the reasoning in the same proof (again applying also to non-orientable maps), trinity symmetry of $M$ implies trinity symmetry of $M^{\sigma}$. Note that both $M$ and $M^{\sigma}$ are non-orientable as their Petrie walks (of length $d$ and $nd$) have odd length. \hfill $\Box$
\medskip

Collecting our findings we arrive at the main result of this paper as a consequence of Theorem \ref{thm:m1} and the remark following it, both in combination with Theorem \ref{thm:acs}.

\begin{theorem}\label{thm:main}
For every odd $d\ge 5$ there exists a regular, self-dual and self-Petrie-dual map of valency $d$. \hfill $\Box$
\end{theorem}

A few remarks are in order. The reader may have observed that if the conclusion of Proposition \ref{prop:pNg} in section \ref{sec:alg} was valid for all odd $k\ge 5$ (and not just for {\em prime} $k\ge 5$), we would have a proof of our main result that would be independent on coverings and the resulting regular maps with trinity symmetry would have automorphism group isomorphic to ${\rm PSL}(2,p)$ for suitable primes depending on $k$. Research in this direction is currently being undertaken by the first two authors of this paper. Here we include a table of the first few values of $N(g)$ for odd $k$ between $5$ and $29$, with $\pi(n)$ standing for the prime factorization of $n$; observe that all the primes $p$ in the prime factorization of $|N(g)|$ satisfy $p\equiv \pm 1\ ({\rm mod}\ 2k\ {\rm and}\ 12)$:

\[\begin{array}{ccccccc}
k &&& N(g) &&& \pi(|N(g)|) \\
\hline
5 &&& -11 &&& {\rm prime} \\
7 &&& -13 &&& {\rm prime} \\
9 &&& -73 &&& {\rm prime} \\
11 &&& +263 &&& {\rm prime} \\
13 &&& -131 &&& {\rm prime} \\
15 &&& -239 &&& {\rm prime} \\
17 &&& -4079 &&& {\rm prime} \\
19 &&& +15503 &&& 37\times 419 \\
21 &&& +5209 &&& {\rm prime} \\
23 &&& -4093 &&& {\rm prime} \\
25 &&& +56149 &&& {\rm prime} \\
27 &&& -16417 &&& {\rm prime} \\
29 &&& +3161869 &&& 59\times 53591 \\
\hline
\end{array}\]
\medskip

\noindent Existence of the corresponding regular maps of valency $k$ with trinity symmetry for odd $k$ between $5$ and $19$ with automorphism group ${\rm PSL}(2,p)$ for the corresponding primes from the above table was computationally verified by M. Conder \cite{Con}, finding such maps of valency $7$ and $17$ also for the Janko simple groups $J_2$ and $J_3$, respectively.
\smallskip

We conclude by noting that a strategy for proving Theorem \ref{thm:main} was also outlined by S. Wilson \cite{W2} by reducing the problem to a number-theoretic question related to Chebyshev polynomials over finite fields.

\bigskip
\bigskip

\noindent{\bf Acknowledgement.}~~ The third author acknowledges support from the APVV Research Grants 15-0220 and 17-0428, and from the VEGA Research Grants 1/0026/16 and 1/0142/17.

\end{document}